\documentclass[12pt, reqno]{amsart}
\usepackage{amsmath, amsthm, amscd, amsfonts, amssymb, graphicx, color, mathrsfs}
\usepackage[bookmarksnumbered, colorlinks, plainpages]{hyperref}
\hypersetup{colorlinks=true,linkcolor=red, anchorcolor=green, citecolor=cyan, urlcolor=red, filecolor=magenta, pdftoolbar=true}
\usepackage{hyperref}
\newtheorem{theorem}{Theorem}[section]
\newtheorem{lemma}[theorem]{Lemma}

\theoremstyle{definition}

\theoremstyle{remark}
\newtheorem{remark}[theorem]{Remark}
\numberwithin{equation}{section}

\begin{document}

\setcounter{page}{1}

\title[Uniqueness for an inverse problem of determining orders and coefficients \dots ]{Uniqueness for an inverse problem of determining orders and coefficients for multi-term time-fractional diffusion equation
}

\author[R. Ashurov and  D. Shamuratov]{Ravshan Ashurov and Damir Shamuratov}

\address{\textcolor[rgb]{0.00,0.00,0.84}{Ravshan Ashurov:
V. I. Romanovskiy Institute of Mathematics Uzbekistan Academy of Sciences,
University street 9, Tashkent--100174,
Uzbekistan; Engineering School, Central Asian University, 111221, Tashkent, Uzbekistan.}}
\email{\textcolor[rgb]{0.00,0.00,0.84}
{ashurovr@gmail.com}}

\address{\textcolor[rgb]{0.00,0.00,0.84}{Damir Shamuratov:
V. I. Romanovskiy Institute of Mathematics Uzbekistan Academy of Sciences,
University street 9, Tashkent--100174,
Uzbekistan}}
\email{\textcolor[rgb]{0.00,0.00,0.84}{damirshamuratov5@gmail.com}}

%\dedicatory{This paper is dedicated to Professor ABCD}

\subjclass[2010]{}


\begin{abstract}
This paper investigates a uniqueness problem for determining all fractional orders and their corresponding coefficients in multi-term time-fractional diffusion equations with Caputo derivatives in an abstract Hilbert space. The main result establishes that these parameters can be uniquely identified from a single scalar observation of the solution, of the form \((u(t),\phi)\), by comparing the solutions of two evolution equations governed by possibly different positive self-adjoint operators \(A\) and \(B\). In contrast to existing results that rely on local space–time
measurements and, in particular, may require the initial data to satisfy sign conditions
and strong spatial regularity assumptions, our approach uses only the
minimal assumption $a \in H$, without imposing any sign or additional smoothness
conditions on the initial value.
The proof combines spectral representations of the solutions with a recursive asymptotic comparison method based on fundamental identities for the multinomial Mittag–Leffler function. This approach enables the successive identification of all fractional orders and their associated coefficients, even when the two equations involve different operators. Furthermore, the uniqueness argument extends to multi-term time-fractional wave equations with orders in \((1,2)\).
\end{abstract}
 \maketitle

\section{Introduction}\label{sec1}

Fractional differential equations have become an important tool for modeling anomalous diffusion, viscoelasticity, heat conduction in heterogeneous media, signal processing, and many other phenomena exhibiting memory and hereditary effects. In recent decades, substantial attention has been devoted to the mathematical analysis of fractional evolution equations, including questions of existence, uniqueness, regularity, asymptotic behavior, controllability, and inverse problems. Among various models, multi-term time-fractional diffusion equations have attracted particular interest because they provide a more accurate description of complex physical processes than their single-term counterparts.

In many practical applications, the parameters appearing in fractional differential equations are not known a priori and must be identified from experimental or observational data. Consequently, inverse problems concerning the determination of fractional orders and associated coefficients have become an active area of research. Such problems are of considerable importance because the fractional orders characterize the memory properties of the underlying medium and therefore play a crucial role in the interpretation of physical processes.

Let $A$ be an arbitrary positive self-adjoint operator defined in a separable Hilbert space $H$ with domain $D(A)$. Denote the inner product in $H$ by $(\cdot,\cdot)$ and the norm by $\|\cdot\|$. Suppose further that the inverse operator $A^{-1}$ is compact. 
Then it is well known that $A$ possesses a complete orthonormal system of eigenfunctions $\{v_k\}$ and a countable set of positive eigenvalues $\{\lambda_k\}$ with a unique accumulation point at $+\infty$. 
We assume that the eigenvalues are ordered in a nondecreasing manner, i.e.,
\[
0 < \lambda_1 \leq \lambda_2 \leq \cdots \to +\infty.
\]

The Caputo fractional derivative of order $0<\alpha_j<1$ $(j=1,\dots,n)$ for a function $h(t)$ is defined as (see \cite{Kilbas}, p. 92):
\[
\partial_t^{\alpha_j} h(t) =
\frac{1}{\Gamma(1-\alpha_j)} 
\int_0^t 
\frac{h'(s)}{(t-s)^{\alpha_j}} ds, \quad t>0,
\]
provided that the right-hand side exists. Here, \(\Gamma(\cdot)\) is Euler’s gamma function.

Consider the following multi-term fractional diffusion equation with Caputo derivatives:
\begin{equation}\label{1.1}
\begin{cases}
\displaystyle
\sum_{j=1}^{n} q_j \partial_t^{\alpha_j} u(t) + Au(t)
=0,\quad\quad 0<t\le T,\\
u(0)=a,
\end{cases}
\end{equation}
where
\[
0<\alpha_n<\alpha_{n-1}<\cdots<\alpha_1<1,
\]
\(q_j>0\) for \(j=1,\dots,n\) (without loss of generality, \(q_1=1\)) and $a\in H$.

For single-term time-fractional diffusion equations, inverse problems for determining the fractional order have been extensively investigated by many authors. Various uniqueness and reconstruction results have been obtained using additional observations at interior points, boundary measurements, integral conditions, and overdetermination data. We refer the reader to \cite{AlimovAshurov1}--\cite{Pskhu1} and the references therein for a comprehensive overview of this topic.

Compared with the single-term case, inverse problems for multi-term time-fractional equations are significantly more challenging because several fractional orders interact simultaneously in the solution representation. As a consequence, the available literature on such problems remains relatively limited.

One of the first systematic studies in this direction was carried out by \cite{LiYamamoto}, where uniqueness results for determining fractional orders in multi-term time-fractional diffusion equations were established by means of eigenfunction expansions and Laplace transform techniques. 

More recently, Jiang and Li \cite{JiangLi} investigated an inverse problem for multi-term time-fractional wave equations and established uniqueness in determining the fractional orders. Their approach is based on local space--time observation data and requires relatively strong assumptions on the initial value. Moreover, their uniqueness result is established under the assumption that two solutions corresponding to different fractional orders coincide on a nonempty space--time subdomain \(\omega\times(0,T)\), where \(\omega\subset\Omega\) is nonempty. More precisely, they consider the following pair of initial-boundary value problems:
$$
\left\{
\begin{aligned}
&\sum_{j=1}^{m} c_j \partial_t^{\alpha_j} u(x,t)
-\Delta u(x,t)
+p(x)u(x,t)
=0,
&& (x,t)\in \Omega\times(0,T),\\
&u(x,0)=a(x),\qquad u_t(x,0)=0,
&& x\in\Omega,\\
&u(x,t)=0,
&& (x,t)\in \partial\Omega\times(0,T),
\end{aligned}
\right.
$$
and
$$
\left\{
\begin{aligned}
&\sum_{j=1}^{m} c_j \partial_t^{\beta_j} v(x,t)
-\Delta v(x,t)
+q(x)v(x,t)
=0,
&& (x,t)\in \Omega\times(0,T),\\
&v(x,0)=a(x),\qquad v_t(x,0)=0,
&& x\in\Omega,\\
&v(x,t)=0,
&& (x,t)\in \partial\Omega\times(0,T),
\end{aligned}
\right.
$$
where $1<\alpha_m<\cdots<\alpha_1<2,$\quad$1<\beta_m<\cdots<\beta_1<2,$ \; $c_j>0,\;(c_1=1),\, j=1,\ldots,m.$ 
The proof of the uniqueness of the potential, namely $p(x)=q(x)$, can be found in \cite{Bazlekova}.
In contrast, the present work adopts a different observation strategy and a more general operator framework. We consider two multi-term time-fractional evolution equations governed by possibly distinct positive self-adjoint operators \(A\) and \(B\) in an abstract Hilbert space \(H\). 

Moreover, to determine the fractional orders and their corresponding coefficients, we consider the following observation data:
\[
(u(t),\phi),\qquad 0<t<T,
\]
for some prescribed weight function \(\phi\). In some existing approaches, uniqueness is established using different types of observation data. For example, in \cite{LiYamamoto}, the observation is given by pointwise data at a fixed spatial point, whereas in \cite{JiangLi}, uniqueness is established by requiring two solutions to coincide throughout a nonempty subdomain \(\omega\subset\Omega\) over the observation interval \((0,T)\). Such spatially distributed pointwise measurements are difficult to obtain in practice, since physical sensors have a finite spatial extent rather than acting as ideal point detectors. Our approach requires only the scalar quantity $(u(t),\phi)$ at each time $t\in(0,T)$ and also allows us to establish uniqueness under weaker regularity assumptions on the initial data $a$. Finally, the weight function $\phi$  may be chosen to vanish outside $\omega$. As is readily seen from the proofs, all the main results remain valid with this choice.

 For the single-term time-fractional diffusion equation, uniqueness results for determining the fractional order from observation data were obtained in \cite{AshurovYamamoto} under such a condition.

The proof is based on the spectral representation of solutions and a detailed asymptotic analysis of the multinomial Mittag--Leffler function. A major difficulty arises from the appearance of multiple competing asymptotic terms generated by different fractional orders. To overcome this difficulty, we repeatedly employ a fundamental identity for the multinomial Mittag--Leffler function and develop a recursive asymptotic comparison procedure. This method allows us to isolate dominant contributions successively and recover each fractional order together with its corresponding coefficient in a step-by-step manner.

The remainder of the paper is organized as follows. In Section \ref{sec2}, we present preliminary results concerning spectral theory and the multinomial Mittag--Leffler function. Section \ref{sec3} is devoted to the formulation and analysis of the inverse problem, where the uniqueness theorem is proved. Finally, concluding remarks are given in Section \ref{sec4}.
\section{Preliminaries}\label{sec2}

Let $\tau$ be an arbitrary real number. We introduce the power of the operator $A$, acting in $H$ according to the rule :
$$
A^{\tau} h = \sum_{k=1}^{\infty} \lambda_k^{\tau} h_k v_k.
$$

Here and everywhere below, for the vector $h \in H$, the symbol $h_k$ denotes the Fourier coefficients of this vector: $h_k = (h, v_k)$. The domain of definition of this operator is determined by the condition $A^{\tau} h \in H$ and has the form:
$$
D(A^{\tau}) = \Big\{ h \in H : \sum_{k=1}^{\infty} \lambda_k^{2{\tau}} |h_k|^2 < \infty \Big\}.
$$

The solution of problem \eqref{1.1} involves the multivariable (multinomial) Mittag--Leffler function. Therefore, we first present the definition of the multinomial Mittag--Leffler function.

The multinomial Mittag--Leffler function is defined as 
\[
E_{(\alpha_1, \dots, \alpha_n),\, \alpha_0}(z_1, \dots, z_n) =
\sum_{k=0}^{\infty} 
\sum_{k_1 + \dots + k_n = k} 
\binom{k}{k_1, \dots, k_n}  
\frac{\prod_{j=1}^n z_j^{k_j}}{\Gamma(\alpha_0 + \sum_{j=1}^n \alpha_j k_j)},
\]
where \(\binom{k}{k_1, \dots, k_n}\) denotes
the multinomial coefficient
\[
\binom{k}{k_1, \dots, k_n} = \frac{k!}{k_1! \cdots k_n!},\quad
\text{with} \quad
k = \sum_{j=1}^{n} k_j,
\]
and \(k_j\), \(1 \le j \le n\), are non-negative integers. 
For more details of the multinomial Mittag-Leffler function, see \cite{MasahiroYamamoto}, \cite{Luchko}.

\begin{lemma}\label{lemma}\cite{Choe}
Let  $\alpha_j >0$ $(j = 1,\ldots,n)$ and $z_j \in \mathbb{C}$ $(j = 1,\ldots,n)$. Then
\[
\frac{1}{\Gamma(\alpha_0)}
+ \sum_{j=1}^{n} z_j \,
E_{(\alpha_1,\ldots,\alpha_n),\,\alpha_0+\alpha_j}(z_1,\ldots,z_n)
=
E_{(\alpha_1,\ldots,\alpha_n),\,\alpha_0}(z_1,\ldots,z_n).
\]
\end{lemma}

For later use, we adopt the abbreviation
$$
E_{\alpha',\beta}(-
\lambda_kt^{\alpha_1},*)=E_{(\alpha_1,\alpha_1-\alpha_2,...,\alpha_1-\alpha_n),\beta}(-\lambda_k t^{\alpha_1},-q_2t^{\alpha_1-\alpha_2},...,-q_nt^{\alpha_1-\alpha_n}),
$$

By Lemma \eqref{lemma}, we obtain the following identity 
\begin{equation}\label{2.1}\begin{aligned}
  E_{\alpha',1+\alpha_1}(-\lambda_kt^{\alpha_1},*)
=\frac{1}{\Gamma(1+\alpha_1)}
&-\lambda_k t^{\alpha_1}
E_{\alpha',1+2\alpha_1}(-\lambda_kt^{\alpha_1},*)\\
&-\sum_{j=2}^{n}
q_j t^{\alpha_1-\alpha_j}
 E_{\alpha',1+2\alpha_1-\alpha_j}(-\lambda_kt^{\alpha_1},*).\end{aligned} 
\end{equation}

\begin{lemma}\label{l2.1} \cite{KarelVanBockstal} Let \(\beta>0\) and \(\alpha_1 > \alpha_2 > \cdots > \alpha_m > 0\) be given with \(\alpha_1 \in (0,2)\). Assume that \(\alpha_1 \pi/2 < \mu < \min\{\alpha_1 \pi, \pi\}\), \(\mu \le |\arg(z_1)| \le \pi\), and there exists \(K>0\) such that \(-K \le z_j < 0\) for \(j=2,\ldots,m\). Then there exists a constant \(C>0\) depending only on \(\mu, K, \alpha_j \ (j=1,\ldots,m)\) and \(\beta\) such that
\[
\left| E_{(\alpha_1,\ldots,\alpha_m),\beta}(z_1,\ldots,z_m) \right|
\le
\frac{C}{1+|z_1|}.
\]
\end{lemma}

\begin{lemma}\cite{MasahiroYamamoto}\label{Solve}
    The solution of \eqref{1.1} is represented by

    $$u(t)=\sum_{k=1}^{\infty}\left(1-\lambda_kt^{\alpha_1}E_{\alpha',1+\alpha_1}(-\lambda_kt^{\alpha_1},*)\right)a_kv_k,
$$
where $(a,v_k)=a_k$, $a\in H.$
\end{lemma}

\section{Inverse Problem}\label{sec3}

Let us formulate our main result. For this, we will look at another system
\begin{equation}\label{3.1}
\begin{cases}
\displaystyle \sum_{j=1}^{n} \mu_j \, \partial_t^{\beta_j} \omega(t) +B\omega(t) = 0,\quad\quad 0<t\le T, \\[1mm]
\omega(0) = a, 
\end{cases}
\end{equation}
where \(0<\beta_n<\beta_{n-1}<\dots<\beta_1< 1\), \(\mu_j>0\) for \(j=1,\dots,n\) (without loss of generality, \(\mu_1=1\)) and $B$ is a positive definite self-adjoint operator with compact inverse. By the spectral theorem, the operator $B$ possesses a complete orthonormal system of eigenfunctions $\{\widetilde{v_k}\}$ corresponding to positive eigenvalues $\{\widetilde{\lambda_k}\}$.
Then there exists a unique solution  \cite{MasahiroYamamoto}:
$$\omega(t)=\sum_{k=1}^{\infty}\left(1-{\widetilde{\lambda_k}}t^{\beta_1}E_{\beta',1+\beta_1}(-{\widetilde{\lambda_k}}t^{\beta_1},*)\right)\widetilde{a_k}\widetilde{v_k}.$$
where $(a,\widetilde{v_k})=\widetilde{a_k}$ and $a\in H.$

\begin{theorem}\label{thm3.1}
    Let  $\phi \in D(A^2)\cap D(B^2)$, and $(a,A\phi)\neq 0$, $(a,B\phi)\neq 0.$       If
\[
(u(t), \phi) = (\omega(t), \phi), \quad 0<t<T,
\] then
\[ (a,A\phi)=(a,B\phi),\quad
\alpha_j = \beta_j, \quad q_j=\mu_j, \quad j=1,\dots,n.
\]

\end{theorem}

\begin{proof}
Using Lemma \ref{Solve} together with formula \eqref{2.1}, we derive the following expression:
\begin{equation}\label{3.2}
\begin{aligned}(u(t), \phi) =\sum_{k=1}^{\infty}a_k(\phi,v_k)-\sum_{k=1}^{\infty}\frac{\lambda_ka_k(\phi,v_k)t^{\alpha_1}}{\Gamma(1+\alpha_1)}+\sum_{k=1}^{\infty}\lambda_k^2a_k(\phi,v_k)t^{2\alpha_1}E_{\alpha',1+2\alpha_1}(-\lambda_kt^{\alpha_1},*) \\
+\sum_{k=1}^{\infty}\lambda_ka_k(\phi,v_k)\sum_{j=2}^{n}q_jt^{2\alpha_1-\alpha_j}E_{\alpha',1+2\alpha_1-\alpha_j}(-\lambda_kt^{\alpha_1},*),\end{aligned}
\end{equation}
and
\begin{equation}\label{3.3}\begin{aligned}(\omega(t), \phi) =\sum_{k=1}^{\infty}\widetilde{a_k}(\phi,\widetilde{v_k})-\sum_{k=1}^{\infty}\frac{\widetilde{\lambda_k}\widetilde{a_k}(\phi,\widetilde{v_k})t^{\beta_1}}{\Gamma(1+\beta_1)}+\sum_{k=1}^{\infty}\widetilde{\lambda_k}^2\widetilde{a_k}(\phi,\widetilde{v_k})t^{2\beta_1}E_{\beta',1+2\beta_1}(-\widetilde{\lambda_k}t^{\beta_1},*) \\
+\sum_{k=1}^{\infty}\widetilde{\lambda_k}\widetilde{a_k}(\phi,\widetilde{v_k})\sum_{j=2}^{n}\mu_jt^{2\beta_1-\beta_j}E_{\beta',1+2\beta_1-\beta_j}(-\widetilde{\lambda_k}t^{\beta_1},*).
\end{aligned}\end{equation}

Since the estimates in \eqref{3.2} and \eqref{3.3} are obtained in the same manner, it suffices to consider only \eqref{3.2} in the subsequent analysis.

We introduce the following notation:
$$I_{1}(t)=\sum_{k=1}^{\infty}\lambda_k^2a_k (\phi,v_k)t^{2\alpha_1}E_{\alpha',1+2\alpha_1}(-\lambda_kt^{\alpha_1},*) ,$$
$$I_{2}(t)=\sum_{k=1}^{\infty}\lambda_ka_k (\phi,v_k)\sum_{j=2}^{n}q_jt^{2\alpha_1-\alpha_j}E_{\alpha',1+2\alpha_1-\alpha_j}(-\lambda_kt^{\alpha_1},*).$$

Using Lemma \ref{l2.1}, we obtain
$$|I_{1}(t)|\le Ct^{2\alpha_1}\sum_{k=1}^{\infty}\lambda_k^2|a_k||(\phi,v_k)|,$$
$$|I_{2}(t)|\le C\sum_{j=2}^{n}q_jt^{2\alpha_1-\alpha_j}\sum_{k=1}^{\infty}\lambda_k|a_k|| (\phi,v_k)|.$$

Using the Cauchy--Schwarz inequality and since $a\in H$,\,$\phi\in D(A^2)\cap D(B^2)$, we have
$$|I_{1}(t)|\le Ct^{2\alpha_1}\left(\sum_{k=1}^{\infty}|a_k|^2\right)^{1/2}\left(\sum_{k=1}^{\infty}\lambda_k^4|(\phi,v_k)|^2\right)^{1/2}= Ct^{2\alpha_1}\|a\|\|A^2\phi\|.$$

Similarly,
$$|I_2(t)|\le C\sum_{j=2}^{n}t^{2\alpha_1-\alpha_j}\left(\sum_{k=1}^{\infty}|a_k|^2\right)^{1/2}\left(\sum_{k=1}^{\infty}\lambda_k^2|(\phi,v_k)|^2\right)^{1/2}=C\sum_{j=2}^{n}t^{2\alpha_1-\alpha_j}\|a\|\|A\phi\|.$$

Therefore,  $\min\limits_{2\le j \le n}\{2\alpha_1-\alpha_j\}=2\alpha_1-\alpha_2$, we get 
$$
(u(t),\phi)=
(a,\phi)-\frac{(a,A\phi)t^{\alpha_1}}{\Gamma(1+\alpha_1)}
+ O\left(t^{2\alpha_1}\right)
+ O\left(t^{2\alpha_1-\alpha_2}\right),\quad t\to 0.
$$
 Given that $a\in H$,\,$\phi\in D(A^2)\cap D(B^2)$, similarly
$$
(\omega(t),\phi)=
(a,\phi)-\frac{(a,B\phi)t^{\beta_1}}{\Gamma(1+\beta_1)}
+ O\left(t^{2\beta_1}\right)
+ O\left(t^{2\beta_1-\beta_2}\right),\quad t\to 0.
$$

Since 
$$(u(t), \phi) = (\omega(t), \phi),$$
we have,
\begin{equation}\label{3.6}
-\frac{(a,A\phi)t^{\alpha_1}}{\Gamma(1+\alpha_1)}+o(t^{\alpha_1})= -\frac{(a,B\phi)t^{\beta_1}}{\Gamma(1+\beta_1)}+o(t^{\beta_1}), \quad t\to 0.
\end{equation}

We assume that $\alpha_1<\beta_1$ (or $\alpha_1>\beta_1$). Dividing both sides of \eqref{3.6} by $t^{\alpha_1}$ (or $t^{\beta_1}$) and passing to the limit as $t\to0$, we arrive at
\[
\frac{(a,A\phi)}{\Gamma(1+\alpha_1)}=0 \left(\text{or}\quad  \frac{(a,B\phi)}{\Gamma(1+\beta_1)}=0\right).
\]

Using the fact that $(a,A\phi)\neq 0$ (or $(a,B\phi)\neq 0$), this is impossible. Hence, the assumption $\alpha_1<\beta_1$ ($\alpha_1>\beta_1$) is false, and consequently
\[
\alpha_1=\beta_1.
\]

It follows from \eqref{3.6} and the equality $\alpha_1=\beta_1$ that
$(a,A\phi)=(a,B\phi).$

In this way, by subtracting the first and second term from both sides of \eqref{3.2},\eqref{3.3}  and repeating the above steps and using Lemma \ref{lemma}, we rewrite
$
E_{\alpha',\,1+2\alpha_1-\alpha_2}
(-\lambda_k t^{\alpha_1},*)
$
in the following form:
$$ E_{\alpha',1+2\alpha_1-\alpha_2}(-\lambda_kt^{\alpha_1},*)
=\frac{1}{\Gamma(1+2\alpha_1-\alpha_2)}
-\lambda_k t^{\alpha_1}
E_{\alpha',1+3\alpha_1-\alpha_2}(-\lambda_kt^{\alpha_1},*)$$
$$
-\sum_{j=2}^{n}
q_j t^{\alpha_1-\alpha_j}
 E_{\alpha',1+3\alpha_1-\alpha_2-\alpha_j}(-\lambda_kt^{\alpha_1},*),$$
and $t\to 0$, we obtain
\begin{equation}\label{second equality}
  \begin{aligned}  &\frac{q_2(a,A\phi)t^{2\alpha_1-\alpha_2}}{\Gamma(1+2\alpha_1-\alpha_2)}+O(t^{2\alpha_1})+O(t^{3\alpha_1-\alpha_2})+O(t^{3\alpha_1-2\alpha_2})+O(t^{2\alpha_1-\alpha_j})\\
    &=\frac{\mu_2(a,B\phi)t^{2\beta_1-\beta_2}}{\Gamma(1+2\beta_1-\beta_2)}+O(t^{2\beta_1})+O(t^{3\beta_1-\beta_2})+O(t^{3\beta_1-2\beta_2})+O(t^{2\beta_1-\beta_j}),\end{aligned}
\end{equation}
(where $j=3,\dots,n$), or
$$\frac{q_2(a,A\phi)t^{2\alpha_1-\alpha_2}}{\Gamma(1+2\alpha_1-\alpha_2)}+o(t^{2\alpha_1-\alpha_2})=\frac{\mu_2(a,B\phi)t^{2\beta_1-\beta_2}}{\Gamma(1+2\beta_1-\beta_2)}+o(t^{2\beta_1-\beta_2}).$$

Having established that $\alpha_1=\beta_1$ and $(a,A\phi)=(a,B\phi)$, we repeat the above procedure. By the same argument, we get
\[
\alpha_2=\beta_2
\quad \text{and} \quad
q_2=\mu_2.
\]

To prove that $\alpha_3=\beta_3$, we subtract from \eqref{3.2} and \eqref{3.3} the first two terms together with
\begin{equation}\label{equal}
\cfrac{q_2(a,A\phi)t^{2\alpha_1-\alpha_2}}{\Gamma(1+2\alpha_1-\alpha_2)}, \qquad
\cfrac{\mu_2(a,B\phi)t^{2\beta_1-\beta_2}}{\Gamma(1+2\beta_1-\beta_2)}.
\end{equation}
Because the first and second terms, also \eqref{equal} terms  ($\alpha_1=\beta_1,\,\alpha_2=\beta_2,\,q_2=\mu_2,\,(a,A\phi)=(a,B\phi)$) are equal. 
Then we obtain
$$
\cfrac{q_3(a,A\phi)t^{2\alpha_1-\alpha_3}}{\Gamma(1+2\alpha_1-\alpha_3)}
+O(t^{2\alpha_1})
+O(t^{3\alpha_1-\alpha_2})
+O(t^{3\alpha_1-2\alpha_2})$$$$
+O(t^{3\alpha_1-\alpha_3})
+O(t^{3\alpha_1-\alpha_3-\alpha_2})+O(t^{2\alpha_1-\alpha_j})
$$
$$
=\cfrac{\mu_3(a,B\phi)t^{2\beta_1-\beta_3}}{\Gamma(1+2\beta_1-\beta_3)}+O(t^{2\beta_1})
+O(t^{3\beta_1-\beta_2})
+O(t^{3\beta_1-2\beta_2})$$$$
+O(t^{3\beta_1-\beta_3})
+O(t^{3\beta_1-\beta_3-\beta_2})+O(t^{2\beta_1-\beta_j}).
$$
where $j=4,\dots,n.$
In order to rewrite this equality in the form
\begin{equation}\label{new}
\frac{q_3(a,A\phi)t^{2\alpha_1-\alpha_3}}{\Gamma(1+2\alpha_1-\alpha_3)}
+o(t^{2\alpha_1-\alpha_3})=
\frac{\mu_3(a,B\phi)t^{2\beta_1-\beta_3}}{\Gamma(1+2\beta_1-\beta_3)}
+o(t^{2\beta_1-\beta_3}),
\quad t\to 0,
\end{equation}
it is required that
$$
3\alpha_1-2\alpha_2>2\alpha_1-\alpha_3,
$$
$$
3\beta_1-2\beta_2>2\beta_1-\beta_3,
$$
which is not always satisfied.
However, by applying Lemma \ref{lemma} $s$ and $p$  times to the multinomial Mittag--Leffler functions
$
E_{\alpha',1+2\alpha_1-\alpha_2}(-\lambda_kt^{\alpha_1},*)
$, $
E_{\beta',1+2\beta_1-\beta_2}(-\lambda_kt^{\beta_1},*)
$
appearing in \eqref{3.2} and \eqref{3.3} there exists $s, p\in \mathbb{N}$ (here
and henceforth let $\mathbb{N}:= \{1, 2, 3, ...\}$) such that
$$
2\alpha_1-\alpha_2+s(\alpha_1-\alpha_2)>2\alpha_1-\alpha_3,
$$
$$
2\beta_1-\beta_2+p(\beta_1-\beta_2)>2\beta_1-\beta_3,
$$
or, more explicitly,
$$
s(\alpha_1-\alpha_2)>\alpha_2-\alpha_3,
$$
$$
p(\beta_1-\beta_2)>\beta_2-\beta_3,
$$
where $s=\min\{\delta\in\mathbb{N}:\delta(\alpha_1-\alpha_2)>\alpha_2-\alpha_3\}$, $p=\min\{\delta\in\mathbb{N}:\delta(\beta_1-\beta_2)>\beta_2-\beta_3\}$.
If $p>s$, for all terms up to $s$, we obtain the following equality when $t\to 0$:
\begin{equation}\label{eq3.8}\begin{aligned}
(a,A\phi)\frac{(-1)^{i}q_2^{i+1}t^{2\alpha_1-\alpha_2+i(\alpha_1-\alpha_2)}}{\Gamma(1+2\alpha_1-\alpha_2+i(\alpha_1-\alpha_2))}
+o(t^{2\alpha_1-\alpha_2+i(\alpha_1-\alpha_2)})\\
=
(a,B\phi)\frac{(-1)^{i}\mu_2^{i+1}t^{2\beta_1-\beta_2+i(\beta_1-\beta_2)}}{\Gamma(1+2\beta_1-\beta_2+i(\beta_1-\beta_2))}
+o(t^{2\beta_1-\beta_2+i(\beta_1-\beta_2)}),\end{aligned}
\end{equation}
where $i=1,2,\dots,s-1$. As shown above, we already know that
$
\alpha_1=\beta_1,\,
\alpha_2=\beta_2,\, q_2=\mu_2,\, (a,A\phi)=(a,B\phi).
$ 

Considering the terms from $s$ to $p-1$, we obtain the following equality:
\begin{equation}\label{eq3.9}
    \begin{aligned}
&(a,A\phi)\frac{q_3t^{2\alpha_1-\alpha_3}}{\Gamma(1+2\alpha_1-\alpha_3)}
+o(t^{2\alpha_1-\alpha_3})\\
&=
(a,B\phi)\sum_{i=s}^{p-1}\frac{(-1)^{i}\mu_2^{i+1}t^{2\beta_1-\beta_2+i(\beta_1-\beta_2)}}{\Gamma(1+2\beta_1-\beta_2+i(\beta_1-\beta_2))}
+o(t^{2\beta_1-\beta_2+i(\beta_1-\beta_2)}),\end{aligned}\end{equation}

In view of the fact that $\alpha_1=\beta_1$, $\alpha_2=\beta_2$, and $p>s$, we have
\[
2\beta_1-\beta_2+i(\beta_1-\beta_2)>2\alpha_1-\alpha_3,
\qquad i=s,s+1,\ldots,p-1.
\]

Therefore, dividing both sides of the \eqref{eq3.9} by $t^{2\alpha_1-\alpha_3}$ and letting $t\to0$, we obtain
\[
\frac{q_3(a,A\phi)}{\Gamma(1+2\alpha_1-\alpha_3)}=0.
\]

 Since $q_3\neq 0$ and $(a,A\phi)\neq 0$, the above equality leads to a contradiction. Therefore, the case $p>s$ is impossible. We next consider the case $s>p$. For all terms up to $p-1$, we get \eqref{eq3.8} equality. From $p$ to $s-1$, we obtain following equality:
\begin{equation}\label{eq3.10}
    \begin{aligned}
&(a,A\phi)\sum_{i=p}^{s-1}\frac{(-1)^{i}q_2^{i+1}t^{2\alpha_1-\alpha_2+i(\alpha_1-\alpha_2)}}{\Gamma(1+2\alpha_1-\alpha_2+i(\alpha_1-\alpha_2))}
+o(t^{2\alpha_1-\alpha_2+i(\alpha_1-\alpha_2)}) \\
&=
  (a,B\phi)\frac{\mu_3t^{2\beta_1-\beta_3}}{\Gamma(1+2\beta_1-\beta_3)}
+o(t^{2\beta_1-\beta_3)}),\end{aligned}\end{equation}
Taking into account that $\alpha_1=\beta_1$, $\alpha_2=\beta_2$, and $s>p$, we have
\[
2\alpha_1-\alpha_2+i(\alpha_1-\alpha_2)>2\beta_1-\beta_3,
\qquad i=p,p+1,\ldots,s-1.
\]

Therefore, dividing both sides of the \eqref{eq3.10} by $t^{2\beta_1-\beta_3}$ and letting $t\to0$, we obtain
\[
\frac{\mu_3(a,B\phi)}{\Gamma(1+2\beta_1-\beta_3)}=0.
\]

A contradiction is obtained in view of $\mu_3\neq 0$ and $(a,B\phi)\neq 0$. Therefore, $s=p$. Due to the fact that $s=p$, all the equalities in \eqref{eq3.8} hold up to $s$ (or $p$).

As shown above, we already know that
$$
\alpha_1=\beta_1,\qquad
\alpha_2=\beta_2,\qquad q_2=\mu_2,\qquad (a,A\phi)=(a,B\phi).
$$
None of these terms contains $\alpha_3$ or $\beta_3$. Therefore, if
\[
3\alpha_1-2\alpha_2>2\alpha_1-\alpha_3,
\]
then \eqref{new} is valid.
On the other hand, if
\[
3\alpha_1-2\alpha_2<2\alpha_1-\alpha_3,
\]
then there exists \(s\,(\text{or}\,\, p)\in\mathbb{N}\), such that 
$$
2\alpha_1-\alpha_2+s(\alpha_1-\alpha_2)>2\alpha_1-\alpha_3,
$$
we may conclude that \eqref{new} is valid.
By employing the same method as used previously, we derive a contradiction, thereby establishing the result that $\alpha_3=\beta_3$ and $q_3=\mu_3.$

 It may happen that there exists an integer $l\in\mathbb{N}$ such that $l=s-1=p-1$, i.e., the $(s-1)$-th (or equivalently $(p-1)$-th) iteration, for which 
\begin{equation}\label{3.8}\begin{aligned}
2\alpha_1 - \alpha_2 + l(\alpha_1 - \alpha_2) = 2\alpha_1 - \alpha_3,\\
2\beta_1 - \beta_2 + l(\beta_1 - \beta_2) = 2\beta_1 - \beta_3,
\end{aligned}
\end{equation}
or, more explicitly,
\[\begin{aligned}
l(\alpha_1 - \alpha_2) = \alpha_2 - \alpha_3,\\
l(\beta_1 - \beta_2) = \beta_2 - \beta_3
\end{aligned}\]
then we obtain the following equality when $t\to 0$:
\begin{equation}\label{3.9}\begin{aligned}
(a,A\phi)t^{2\alpha_1 - \alpha_2 + l(\alpha_1 - \alpha_2)}\left(\frac{q_3}{\Gamma(1+2\alpha_1-\alpha_3)}+\frac{(-1)^lq_2^{l+1}}{\Gamma(1+2\alpha_1-\alpha_2+l(\alpha_1-\alpha_2))}\right)\\
+o(t^{2\alpha_1 - \alpha_2 + l(\alpha_1 - \alpha_2)})\\=(a,B\phi)t^{2\beta_1 - \beta_2 + l(\beta_1 - \beta_2)}\left(\frac{\mu_3}{\Gamma(1+2\beta_1-\beta_3)}+\frac{(-1)^{l}\mu_2^{l+1}}{\Gamma(1+2\beta_1-\beta_2+l(\beta_1-\beta_2))}\right)\\
+o(t^{2\beta_1 - \beta_2 + l(\beta_1 - \beta_2)}).\end{aligned}\end{equation}

Due to the fact that $\alpha_1=\beta_1$ and $\alpha_2=\beta_2$, it follows from \eqref{3.8} that $$\alpha_3=\alpha_2-l(\alpha_1-\alpha_2)=\beta_2-l(\beta_1-\beta_2)=\beta_3.$$

Moreover, taking into account that $\alpha_1=\beta_1$ , $\alpha_2=\beta_2$, $\alpha_3=\beta_3$ and $(a,A\phi)=(a,B\phi)$ we divide both sides of \eqref{3.9} by
$
t^{2\alpha_1-\alpha_2+l(\alpha_1-\alpha_2)}
$ (or $
t^{2\beta_1-\beta_2+l(\beta_1-\beta_2)}
$)
and letting $t\to 0$, we obtain
\begin{equation}\label{3.10}
    \frac{q_3-\mu_3}{\Gamma(1+2\alpha_1-\alpha_3)}+\frac{(-1)^l(q_2^{l+1}-\mu_2^{l+1})}{\Gamma(1+2\alpha_1-\alpha_2+l(\alpha_1-\alpha_2))}=0.
\end{equation}
We already know that $q_2=\mu_2$, it follows from \eqref{3.10} that $q_3=\mu_3$.

To prove that $\alpha_4=\beta_4$, we subtract from \eqref{3.2} and \eqref{3.3} the first two terms and \eqref{equal} together with
\begin{equation}
\cfrac{q_3(a,A\phi)t^{2\alpha_1-\alpha_3}}{\Gamma(1+2\alpha_1-\alpha_3)}, \qquad
\cfrac{\mu_3(a,B\phi)t^{2\beta_1-\beta_3}}{\Gamma(1+2\beta_1-\beta_3)}.
\end{equation}

Then we obtain
$$\begin{aligned}
\cfrac{q_4(a,A\phi)t^{2\alpha_1-\alpha_4}}{\Gamma(1+2\alpha_1-\alpha_4)}
&+O(t^{2\alpha_1})
+O(t^{3\alpha_1-\alpha_2})
+O(t^{3\alpha_1-2\alpha_2})+O(t^{2\alpha_1-\alpha_j})
\\&+O(t^{3\alpha_1-\alpha_3})
+O(t^{3\alpha_1-\alpha_3-\alpha_2})
+O(t^{3\alpha_1-\alpha_4})+O(t^{3\alpha_1-\alpha_2-\alpha_4})\end{aligned}$$
$$
\begin{aligned}
=\cfrac{\mu_4(a,B\phi)t^{2\beta_1-\beta_4}}{\Gamma(1+2\beta_1-\beta_4)}&+O(t^{2\beta_1})
+O(t^{3\beta_1-\beta_2})
+O(t^{3\beta_1-2\beta_2})+O(t^{2\beta_1-\beta_j})
\\&+O(t^{3\beta_1-\beta_3})
+O(t^{3\beta_1-\beta_3-\beta_2})
+O(t^{3\beta_1-\beta_4})+O(t^{3\beta_1-\beta_2-\beta_4}).
\end{aligned}$$
where $j=5,\dots,n.$
In order to rewrite this equality in the form
\begin{equation}\label{3.12}
\frac{q_4(a,A\phi)t^{2\alpha_1-\alpha_4}}{\Gamma(1+2\alpha_1-\alpha_4)}
+o(t^{2\alpha_1-\alpha_4})=
\frac{\mu_4(a,B\phi)t^{2\beta_1-\beta_4}}{\Gamma(1+2\beta_1-\beta_4)}
+o(t^{2\beta_1-\beta_4}),
\quad t\to 0,
\end{equation}
it is required that
$$
3\alpha_1-\alpha_2-\alpha_3>3\alpha_1-2\alpha_2>2\alpha_1-\alpha_4,
$$
which is not always satisfied. However, similarly to the above case, there exists $m \in \mathbb{N}$ such that
\[
2\alpha_1 - \alpha_2 + m(\alpha_1 - \alpha_2) > 2\alpha_1 - \alpha_4.
\]

Nevertheless, when repeating the procedure for \(m\) times, there may exists \(r\in\mathbb{N}\) such that
\[
2\alpha_1-\alpha_2+r(\alpha_1-\alpha_2)>3\alpha_1-\alpha_3-\alpha_2.
\]

In this case, by applying Lemma \ref{lemma} to \(E_{\alpha',1+2\alpha_1-\alpha_3}(-\lambda_kt^{\alpha_1},*)\) exactly \(h\) times, we can find an integer \(h\) satisfying
\[
2\alpha_1-\alpha_3+h(\alpha_1-\alpha_2)>2\alpha_1-\alpha_4.
\]

Throughout the construction of $m$ and $h$ no terms involving $\alpha_4$ or $\beta_4$ arise, similarly to the situation in \eqref{eq3.8}.
So we may conclude that \eqref{3.12} is valid. By employing the same method as used previously, we derive a contradiction, thereby establishing the result that
$
\alpha_4=\beta_4,\,\,$ 
$
q_4=\mu_4.
$

The subsequent steps follow exactly the same procedure as that used in the proof of
$
\alpha_5=\beta_5
$
and
$
q_5=\mu_5.
$
Repeating the above argument recursively, we conclude that
\[
\alpha_j=\beta_j,
\qquad
q_j=\mu_j,
\qquad
j=1,\ldots,n.
\]

Thus the proof of Theorem \ref{thm3.1} is complete.\end{proof}

\begin{remark}
By reviewing the proof, one can readily verify that the same arguments remain valid for multi-term time-fractional wave equations. More precisely, the results continue to hold for equations \eqref{1.1} and \eqref{3.1} with
\[
1<\alpha_n<\alpha_{n-1}<\dots<\alpha_1<2,
\qquad
1<\beta_n<\beta_{n-1}<\dots<\beta_1<2,
\]
supplemented with the initial conditions
\[
u(0)=a,\,\,\,\,\partial_t u(0)=0  ,
\qquad
w(0)=a,\,\,\,\,\partial_t \omega(0)=0 .
\]
\end{remark}

\section{Conclusion}\label{sec4}

In this paper, we investigated an inverse problem for identifying all fractional orders and their corresponding coefficients in a multi-term time-fractional diffusion equation with Caputo derivatives. The main result establishes the uniqueness of recovering all fractional parameters from a single scalar observation of the solution, namely $(u(t),\phi)$, in an abstract Hilbert space framework.
The proof relies on the spectral representation of the solution together with a careful asymptotic analysis of the multinomial Mittag--Leffler function. A key difficulty stems from the interaction of several asymptotic terms generated by different fractional orders. This difficulty is overcome by repeatedly applying an identity for the multinomial Mittag--Leffler function and developing a recursive asymptotic comparison argument, which makes it possible to separate the dominant contributions and identify the unknown parameters successively.
Compared with the existing literature, the proposed approach significantly weakens the assumptions required for uniqueness. In particular, the initial value is required only to belong to the underlying Hilbert space, without any sign or additional regularity assumptions. The observation function, however, must satisfy suitable regularity and nondegeneracy conditions. These features substantially broaden the applicability of the obtained uniqueness theorem.
Finally, we can see that the same recursive argument extends naturally to multi-term time-fractional wave equations with fractional orders in $(1,2)$.

\section*{Acknowledgments}

The authors are grateful to professors M. Yamamoto and  Z. Sobirov for discussions of these results.

\end{document}